\magnification 1200
\def \sn{{\smallskip \noindent}}
\def \bn {{\bigskip \noindent}}

\def \la {{\longrightarrow}}
\def \O{{\cal O}}

\def \I{{\cal I}}
\def \P{{\bf  P}}
\def \Q{{\bf Q}}
\def \bP{{\bf P}}

\def \E{{\cal E}}
\def \bR{{\bf R}}

\def \C{{\bf C}}

\def \Z{{\bf Z}}

\def \bE{{\bf E}}

\def \L{{\cal L}}
\def \V{{\cal V}}

\font\medium=cmbx10 scaled \magstep1
\font\large=cmbx10 scaled \magstep2

{\large \centerline {Projective Manifolds with splitting Tangent Bundle,\ I}}
\vskip.5cm

{\medium \centerline {Fr\'ed\'eric Campana and Thomas Peternell}}
\vskip.5cm
\vskip.5cm
{\medium \centerline {Introduction}}

\bn A fundamental theorem in differential geometry due to de Rham investigates complete Riemannian manifolds $X$ whose 
tangent bundles split into a direct sum of subbundles which all are invariant under the holonomy. Under these assumptions the universal
cover of $X$ splits as a product of Riemannian manifolds and this splitting is compatible with the splitting of
the tangent bundle of $X.$ The same holds in the complex case for K\"ahler manifolds.  

\bn In this paper we investigate the more general situation when $X$ is a complex projective manifold whose tangent
bundle splits: $ T_X = \bigoplus E_i$. However we do not make any assumption on holonomy or on integrability
of the subbundles $E_i.$ 

\bn {\bf Problem} Does the universal cover $\tilde X$ of $X$ split, i.e. is there a decomposition $\tilde X = \prod_{i=1}^k A_i$
with complex manifolds $A_i$?

\bn Note that we should not ask for a relation between the splitting of $\tilde X$ and the splitting of $T_X$ because
the $E_i$ might not be integrable. If however the $E_i$ are integrable, then we can ask additionally whether the 
decomposition of $T_{\tilde X}$ is induced by the decomposition of $T_X.$ In that case we speak of a {\it diagonal splitting}
(with respect to the splitting on $T_X).$

\bn Of course this problem makes sense also for compact K\"ahler manifolds  but our
results will mainly concern the projective case. Surprisingly, this very fundamental problem has not yet been studied much.

Beauville [Be99] proved that the problem has a positive solution if $X$ is K\"ahler-Einstein or if $\dim X = 2$  ($X$ K\"ahler).
He also pointed out that in general, with $X$ just being a compact manifold, the problem will have a negative
answer. Hopf manifolds provide the easiest examples. 

Next, Druel [Dr99] investigated projective manifolds whose tangent bundles decompose into a direct sum of line 
bundles $L_i$. He gave a positive solution of the problem, even in the strong form, if either all  subbundles formed
by direct sums of the $L_i$ are integrable or if $X$ is minimal, i.e. $K_X$ is nef (e.g. $X$ has seminegative 
Ricci curvature).

\bn The main results of this paper ar as follows. 

\bn {\bf Proposition} {\it Let $X$ be a smooth projective threefold with a splitting $T_X = L \oplus V$ where $L$ is a line bundle.
Let $\varphi: X \la Y$ be a birational contraction in the sense of Mori theory, i.e. $-K_X$ is $\varphi-$ample. 
Then $Y$ is smooth, $\varphi$ is the blow-up along a smooth curve; $L' = \varphi_*(L) $ 
and $V' = \varphi_*(V)^{**}$
are locally free with $L = \varphi^*(L')$ such that
$$ T_Y = L' \oplus V'.$$
Moreover if the universal cover $\tilde Y $ of $Y$ splits (diagonally with respect to $T_Y = L' \oplus V'$), then
$\tilde X$ splits (diagonally with respect to $T_X = L \oplus V).$ }

The content of this proposition is that (by Mori theory) either we are reduced to uniruled threefolds admitting a Mori contraction
to a curve or a surface or that $K_X$ is nef. In the first case we have

\bn {\bf Theorem}  {\it Let $X$ be a smooth projective threefold with $\kappa (X) = - \infty.$ Suppose $T_X = L
\oplus V.$ Then $\tilde X$ splits. The splitting is diagonal with respect to $T_X = L \oplus V$ unless the
following holds: $X$ is the successive blow-up along smooth curves in a smooth projective threefold $Y$,
the splitting $T_X = L \oplus V$ induces canonically a splitting $T_Y = L' \oplus V'$ and there is a 
$\P_1-$bundle structure $\psi: Y \la Z$ such that $L' = T_{Y/Z}.$ }

\bn We also have general results on projective bundles in any dimension, see sect.2.

\bn In case $K_X$ nef, the Beauville decomposition theorem takes care of the case $c_1(X) = 0.$ Here we concentrate on
the case $\dim X = 3$ and $\kappa (X) = 3,$ so that $K_X$ is big and nef. The intermediate case $\kappa (X) = 1,2$ will be treated in the
second part to this paper. If now $K_X$ is ample, then $X$ admits a K\"ahler-Einstein metric, and Beauville's result
in [Be99] settles the problem. So suppose that $K_X$ is not ample; we then go to the canonical model and analyse this
singular variety to obtain

\bn {\bf Theorem} {\it Let $X$ be a smooth projective threefold with $K_X$ big and nef; let $\phi: X \la Y$ be the canonical model.
{\item (1) The universal cover $\tilde Y$ of $Y$ is of the form $\tilde Y \simeq \Delta \times S$ with $\Delta \subset \C$ 
the unit disc and $S$ a surface with only rational double points as singularities.
\item (2) If $g: \hat S \la S$ denotes the minimal resolution, then the universal cover $\tilde X$ of $X$ is of the form 
$$\tilde X \simeq \Delta \times
\hat S$$ and $\phi \simeq {\rm id}_{\Delta}Ê\times g.$  This decomposition is compatible with the decomposition
$T_X = L \oplus V,$ i.e. $\tilde X$ splits diagonally with respect to $T_X = L \oplus V.$ }}

\bn In higher dimensions we study Fano manifolds. Since these are simply connected, we only need to study a decomposition
$T_X = E_1 \oplus E_2$ into 2 factors. Moreover it is easy to see that $b_2(X) \geq 2,$ so that Mori theory can be used 
again. We obtain:

\bn {\bf Theorem} \ Ê{\it Let $X$ be a Fano $n-$fold. Assume that $n \leq 5$ or that every contraction of an extremal ray contracts
a standard rational curve (see Def. 1.1). If $T_X = E_1 \oplus E_2,$
then $X \simeq Z_1 \times Z_2$ diagonally. }

\bn There is an interesting problem on rational curves which is the obstacle to extend this theorem to higher dimensions.
We say that a rational curve $f: \bP_1 \la X$ is standard if 
$$f^*(T_X) = \O(2) \oplus \bigoplus \O(a_i) $$
with all $a_i \leq 1$ (but $a_i$ might be negative). 
It is well known that uniruled varieties always contain standard rational curves, even with all $a_i \geq 0.$ 
Now the problem reads:

\bn {\bf Problem}  Let $X$ be a projective manifold, $\varphi: X \la Y$ an extremal contraction, i.e. $-K_X$ is
$\varphi-$ample. Is there a
standard rational curve contracted by $\varphi?$

If $\dim Y < \dim X,$ then the answer is yes by the last remark. The answer is also positive if $\varphi$ is the
blow-up of a submanifold, but in general the fibers of a birational contraction can be very bad, so that 
new arguments are needed. If this problem has a positive solution, then the last theorem holds in any dimension.

\bn \bn

{\medium \centerline {0. Notations and basic material}}

\bn {\bf 0.1 Notations}
{\item {(1)} The universal cover of a complex manifold $X$ will be always be denoted by $h: \tilde X \la X.$ 
\item {(2)} Given a complex manifold $X$ with a splitting $T_X = \bigoplus_{i=1}^k E_i,$ we let $r_i$ be the rank 
of the vector bundle $E_i$ and denote $L_i = \det E_i.$ 
\item {(3)} A vector bundle $E$ on a projective manifold is {\it almost nef}, if there is a countable family $(A_i)$ of proper
subvarieties $A_i \subset X$, such that $E \vert C$ is nef for all curves $C \not \subset \bigcup_i A_i.$}

\bn {\bf 0.2 Definition}  Let $X$ be a compact complex manifold.
{\item {(1)} We say that $\tilde X$ {\it splits}, if there is a decomposition $\tilde X = \prod_{i=1}^k A_i$ with 
complex manifolds $A_i$.
\item {(2)} We say that $\tilde X$ {it splits diagonally}, if there {\it exists} a decomposition $\tilde X = \prod A_i$
and a decomposition $ T_X = \bigoplus E_i$ such that $h^*(E_i) = p_i^*(T_{A_i})$, where $p_i$ is the projection
onto $A_i$. We will often write $T_{A_i}$ instead of $p_i^*(T_{A_i}).$ 
 
\item {(3)}  Suppose that $T_X = \bigoplus E_i.$ We say that $\tilde X$ {\it splits diagonally with respect to $T_X =
\bigoplus E_i$,} if $\tilde X = \prod A_i$ such that $T_{A_i}Ê= h^*(E_i).$ }

The following observation is trivial but important

\bn {\bf 0.3 Obvious Fact} Suppose $T_X = \bigoplus E_i$. If $\tilde X$ splits diagonally w.r.t the splitting of $T_X,$ then
all bundles $\bigoplus_{i \in I} E_i$ are integrable ($I \subset  \{1, \ldots,k\}).$  

\bn Given a splitting $T_X = \bigoplus_{i=1}^k E_i,$ we obtain a direct sum decomposition
$$ \Omega^q_X = \bigoplus_{i_1 + \ldots + i_k = q} \bigwedge^{i_1}E_1^* \otimes \ldots \otimes \bigwedge^{i_k}E_k^*.$$
Hence $\bigwedge_kE_j^*$ is a subbundle of $\Omega^q_X.$ Then we have

\bn {\bf 0.4 Lemma} {\it Let $X$ be a compact K\"ahler manifold with $T_X = \bigoplus E_i.$ Then 
$$ c_m(E_i) \in H^m(X,\bigwedge^m E_i^*) \subset H^m(X,\Omega^m_X)$$ for all $m$ and all $i$. }

\bn For the {\bf Proof} we essentially refer to [Be99], Lemma 3.1. The reason is that $c_m(E_i)$ can be computed from
the Atiyah class $at(E)$ and that $at(E)$ is in $H^1(X,E_i^* \otimes {\rm End}(E_i)).$  

\bn {\bf 0.5 Theorem} {\it Let $X$ be a compact K\"ahler manifold, $T_X = \bigoplus E_i$ a splitting. Then $\tilde X$ splits
diagonally w.r.t. to the given splitting of $T_X$ if one of the following holds.
{\item {(1)} (Beauville) $X$ is K\"ahler-Einstein;
\item (2) (Beauville) $\dim X = 2;$
\item (3) (Druel) $X$ is projective, all $E_i$ have rank 1 and one of the following holds: $K_X$ is nef or 
$\bigoplus_{i \in I}E_i$ is integrable for all $I \subset  \{1,2, \ldots ,k\}.$ }}

\bn The proof are found in [Be99] and [Dr99].

\bn {\bf 0.6 Notation} {{\item {(1)} An {\it orbifold}Ê(of dimension $n$) is a normal complex space  $X$,
such that every $x \in X$ has an open neighborhood $U$ and a biholomorphic map $h_U : U \la V/G$ where $V \in \C^n$ is an open
set and $G \subset Gl(n,\C)) $ is a finite group such $U$ is $G-$stable.}
{\item {(2)}ÊThe tangent bundle (sheaf $T_X$ of $X$ is by definition the dual of the sheaf of K\"ahler differentials;
actually we have $T_X = i_*(T_{X_{\rm reg}})$ where $i : X_{\rm reg}Ê\la X$ is the
inclusion map. }
{\item {(3)}ÊA K\"ahler orbifold metric on $X$ is a K\"ahler metric $g$ on $X_{\rm reg}$ with K\"ahler form $\omega$ such that
$h^*(\omega)$ extends to $V$ for all charts $h$ as in (1).}}

\bn Then we have the following theorem of de Rham:

\bn {\bf 0.7 Theorem}Ê {\it Let $X$ be a simply connected complete K\"ahler orbifold. Then 
$$ X \simeq \prod X _i$$
where the $X_i$ are irreducible complete K\"ahler orbifolds, i.e. they have irreducible holonomy. }

\bn The proof is just a tedious translation of the arguments given in [KN63] to the orbifold case. We leave
any details to the reader.

\bn \bn

{\medium \centerline {1. Birational contractions and standard rational curves}}

\bn In this section we study projective manifolds with splitting tangent bundle which admit a birational
Mori contraction. We fix the splitting 
$$ T_X = \bigoplus _{i=1}^k E_i.$$
\sn The following notion of standard rational curves will play an important role in our considerations.

\bn {\bf 1.1 Definition } A rational curve $C \subset X$ with normalisation $f: \P_1 \la X$ is {\it standard}Ê if 
$$f^*(T_X) = \O(2) \oplus \bigoplus \O(a_i) $$
with all $a_i \leq 1$ (but $a_i$ might be negative). 

\bn {\bf 1.2 Lemma}  {\it Let $X$ be a smooth uniruled variety. Then there exists a standard rational curve $C \subset X$
with $T_X \vert C$ nef.}

\bn {\bf Proof.} Hwang-Mok [HM98] and Koll\'ar [Ko96,IV.2.10].

\bn One of the special features of standard rational curves in the context of splitting tangent bundles is  given by

\bn {\bf 1.3 Lemma}  {\it Suppose $T_X = \bigoplus E_i.$ Let $C \subset X$ be a standard rational curve. 
After possible renumeration we may assume $f^*(E_j) = \bigoplus \O(a_i)$  with all $a_i \leq 1$  for $j \geq 2$, 
hence $T_C \subset E_1 \vert C.$
Then $L_j \cdot C = 0,$ $(j < k)$, where $L_i = \det E_i.$ 
\sn If $T_X \vert C$ is even nef, i.e. all $a_i \geq 0,$ then
$E_j \vert C = \O_C^{\oplus r_1}$ for $k > j \geq 2.$  }

\bn Here we started the abuse of notation to omit the normalisation $f^*.$ 

\bn {\bf Proof.} Since all $a_i \leq 1,$ we have $H^1(E_j^*\vert C) = 0,$ hence $c_1(E_j^*\vert C) = 0 $ by (0.4),
i.e. $L_j \cdot C = 0$ for $j < k.$ The rest is clear. 

\bn The following is an easy consequence of (1.2). Here we use the notion of an almost nef vector bundle in the sense of
[DPS00]: $E$ is almost ample if there exists a countable union $B$ of proper subvarieties such that $E \vert C$ is nef
for every curve $C$ not contained in $B.$  

\bn {\bf 1.4 Proposition} {\it Let $X$ be a projective manifold, $\varphi: X \la Y$ the contraction of an extremal ray (or face).
Suppose that $\varphi$ has (ideal-theoretically) a non-trivial smooth fiber $F$ or set-theoretically a smooth fiber
$F$ with almost nef normal bundle. Then there exists a standard rational curve $C \subset X$ with
$\dim \varphi (C) = 0.$ } 

\bn {\bf Proof.}  Since $-K_X $ is $\varphi-$ample, $F$ is Fano. Take a standard rational curve $C \subset F $ with $T_F \vert C$
nef. Consider the exact sequence
$$ 0 \la T_F \vert C \la T_X \vert C \la N_{F/X} \vert C  \la 0.  \eqno (*) $$
By our assumptions on $F,$ the conormal bundle $N^*_F \vert C$ is either spanned by global sections or almost nef,
hence $N_{F/X}^*\vert C$ is nef  
so that (*) splits and our claim follows (one can also argue with any extremal rational curve).

\bn {\bf 1.5 Corollary} Ê{\it Let $X$ be a smooth projective threefold and $\varphi: X \la Y$ an extremal contraction.
Then there exists a standard rational curve $C$ contracted by $\varphi.$ }

\bn {\bf Proof.} The proof goes either by direct inspection of Mori's classification of threefold contractions [Mo82]
or by applying (1.4) and by observing that the only case where $\varphi$ does not have set-theoretically a smooth fiber
with the required normal bundle (for applying (1.4))  is when
$\varphi $ contracts a quadric cone. In this last case we can just take $C$ to be a line in the quadric cone. 

\bn {\bf 1.6 Problem}  Let $X$ be a projective manifold, $\varphi: X \la Y$ an extremal contraction. Is there a
standard rational curve contracted by $\varphi?$

\bn {\bf (1.7)} A possible approach to (1.6) is as follows. 
We need to find $V \subset  {\rm Hom}(\P_1,X) $ generically unsplit (in the sense of Koll\'ar [Ko96]) such that the curves of
the family $(C_t)$ defined by $V$ are contracted to points by $\varphi $ and such that
$$ \dim V \geq \sum_{i=1}^k a_i + \dim {\rm Locus}(V).  \eqno (*) $$ 
Here the $a_i$ are just the positive (or non-negative) entries of $f_t^*(T_X),$ where $f_t$ is the normalisation
of $C_t$ and $t$ is general. Moreover ${\rm Locus}(V)$ denotes the subvariety filled up by the family $(C_t).$ 
\sn Observe that if we take just some 
family satisfying (*), then by [Ko96,IV.2.4] we can pass to 
a generically unsplit family, but then possibly (*) is violated for this new family. 
\sn Notice that if $F$ is smooth, this can be achieved just  by taking a (generically) free family whose general
member is standard ([Ko96,IV.2.9]).
\sn Now suppose that (*) holds. First we observe that 
$$ \dim {\rm Locus}(V,0,x)Ê\leq {\rm card}\{i \vert  a_i \geq 1\}.  \eqno (**)$$
Here ${\rm Locus}(V,0,x)$ denotes the set of $f \in V$ with $f(0) = x \in X.$ 
Moreover by [Ko96,IV.2.5]: 
$$ \dim V = \dim {\rm Locus}(V) + \dim {\rm Locus}(V,0,x) + 1. \eqno (***)$$
Then (**) and (***) yield
$$ \dim V \leq \dim {\rm Locus}(V) + {\rm card}\{i \vert a_i \geq 1\} + 1,$$
hence (*) gives 
$$ {\rm card}\{i \vert a_i \geq 1\} + 1 \geq \sum_{i=1}^k a_i.$$
This in turn implies say $a_1 = 2$ and $a_i \leq 1 $ for $i \geq 2.$ 

\bn {\bf 1.8 Proposition} {\it Let $X$ be a projective manifold with $T_X = \bigoplus E_i.$ Let $\varphi: X \la Y$
be an extremal contraction. 
{\item {(1)} If $L_j = {\rm det}E_j$ is not $\varphi-$trivial, then $\dim F \leq r_j = {\rm rk}E_j$
for all fibers $F$ of $\varphi.$ 
\item {(2)} $ \dim F \leq {\rm max}(r_1, \ldots, r_k).$
\item {(3)} Suppose some fiber of $\varphi $ contains a standard rational curve (e.g. $\dim Y < \dim X$ or 
$\varphi$ has a smooth \ non-trivial fiber). Then all but one of the $L_i$ are $\varphi-$trivial, and the remaining one is $\varphi-$ample.}}

\bn {\bf Proof.} (1) If $\dim F > r_j,$ then, noticing that $L_j \vert F $ is either ample or its dual is ample, we obtain
$$ c_1(L_j \vert F)^{r_j+1} \ne 0.$$
This contradicts (0.4). 
\sn (2) This is clear: $-K_X = L_1 \otimes \ldots \otimes L_k,$ on the other hand $-K_X$ is $\varphi-$ample, so not all $L_i$ can be
$\varphi-$trivial. 
\sn (3) This is again clear by (1.3).

\bn {\bf 1.9 Corollary} \ Ê{\it Let $X$ be a smooth projective threefold with a splitting $T_X = L \oplus V.$ Let $\varphi: X \la Y$
be a birational extremal contraction with exceptional divisor $D.$ Then $\dim \varphi(D) = 1,$ so that $Y$ is smooth and
$\varphi$ is the blow-up along the smooth curve $\varphi(D).$ }

\bn {\bf Proof.} Assume $\dim \varphi(D) = 0.$ Then by (1.8), $L$ cannot be $\varphi-$trivial. Let $l$ be a standard rational
curve contracted by $\varphi$ (1.5). Then $T_X \vert l = \O(2) \oplus \O(a) \oplus \O(b)$ with $a,b \leq 1.$ 
Then $L \vert l = T_l$ by (0.4). Now $D$ is either $\P_2$, a smooth quadric or a quadric cone. In the first two cases
we take $l$ to be a line or a ruling line and immediately get a contradiction (there are lines through one point,
having different tangent directions). In the quadric cone case we take $l$ to be a line through the vertex and perform
the same argument.

\bn {\bf 1.10 Proposition} {\it Let $X$ be a smooth projective threefold with a splitting $T_X = L \oplus V.$ 
Let $\varphi: X \la Y$ be a birational contraction. Then $Y$ is smooth, $\varphi$ is the blow-up along a smooth curve; $L' = \varphi_*(L) $ 
and $V' = \varphi_*(V)^{**}$
are locally free with $L = \varphi^*(L')$ such that
$$ T_Y = L' \oplus V'.$$
Moreover if the universal cover $\tilde Y $ of $Y$ splits (diagonally with respect to $T_Y = L' \oplus V'$), then
$\tilde X$ splits (diagonally with respect to $T_X = L \oplus V).$ }

\bn {\bf Proof.} (1) By (1.9) we know that $Y$ is smooth and that $\varphi$ is the blow-up along the smooth curve $\varphi(D).$
Let $l$ be a non-trivial fiber of $\varphi.$ Then $T_X \vert l = \O(2) \oplus \O \oplus \O(-1)$. Now (1.3) implies that
$L \vert l = \O$ and $T_C \subset V \vert C.$ In particular $L = \varphi^*(L')$ with some line bundle $L'$ on $Y$
(necessarily $L' = \varphi_*(L)$). It is moreover clear that $L' \subset T_Y$ is a subbundle ($L$ has no tangent direction
on the fibers of $\varphi.$) We also have
$$ L' \oplus \varphi_*(V) = \varphi_*(L) \oplus \varphi_*(V) = \varphi_*(T_X) \subset T_Y.$$
Thus
$$ T_Y = (\varphi_*(T_X))^{**} = L' \oplus (\varphi_*(V))^{**} = L' \oplus V'. $$
In particular $V'$ is locally free. Now the splitting maps $T_X \la L$ and $T_X \la V$ yield maps 
$T_Y \la L'$ and $T_Y \la V'$ which give a splitting first on $Y \setminus \varphi(D),$ hence everywhere.
\sn (2) It remains to prove the assertion on the universal cover. First we do the diagonal case. So assume $\tilde Y = A' \times B'.$ Since
$\pi_1(X) = \pi_1(Y), $ we have $\tilde X = X \times_Y \tilde Y.$ 
Suppose $\tilde Y = A' \times B'$ and, letting $g: \tilde Y \la Y$ be the projection,  we moreover have $T_{A'} = g^*(L')$
and $T_{B'}Ê= g^*(V').$ Let
$$ \tilde X = X \times_Y \tilde Y$$
with projection $f: \tilde X \la X$, then $\tilde X$ is the universal cover of $X$ and the induced map $\tilde \varphi: \tilde X
\la \tilde Y$ is the blow-up of $Z = g^{-1}(\varphi (D)).$ 
Then $\varphi^{-1}(A' \times \{b\})$ and $\varphi^{-1}(\{a\}Ê\times B')$ are the leaves of $f^*(L)$ resp. $f^*(V)$. This
gives clearly a decomposition $\tilde X = A \times B$ with $A = A'$ and with $B = \{a\}Ê\times B \la \{a\} \times B' = B'$
the blow-up of $Z \cap \{a\} \times B'.$ 
\sn (3) Finally we have to do the case that the splitting $\tilde Y = A' \times B'$ is not diagonal, i.e. $V$ is not integrable.
(to be done).

\bn {\bf (1.11)} As conclusion we have reduced the study of threefolds $X$ with $T_X = L \oplus V$ to the two following cases:
{\item {(a)} $K_X$ is nef;
\item {(b)} $X$ admits a contraction $\varphi: X \la Y$ with $\dim Y < \dim X;$ in particular $X$ is uniruled.}

\bn \bn {\medium \centerline {2. Projective bundles, conic bundles and del Pezzo fibrations}}

\bn We consider again a smooth projective threefold $X$ with $T_X = L \oplus V$. We suppose moreover that there is
a contraction $\varphi: X \la Y$ with $\dim Y = 2$ or $\dim Y = 1,$ i.e. we have a fiber space by conics or del Pezzo
surfaces. We begin with the conic bundle case. We shall see that the contractions are all submersion, therefore we
are lead to consider projective bundles. We shall start with these and will not need any dimension restriction  here. 

\bn {\bf 2.1 Definition} A vector bundle $E$ of rank $r$ on the projective (or K\"ahler) manifold $X$ is called {\it numerically flat},
if the $\Q-$bundle $E \otimes {{\det E^*} \over {r}}$ is nef, which is the same as to say that $S^rE \otimes \det E^*$ is nef.

\bn For details on numerically trivial vector bundles we refer to [DPS94], in particular for a proof of

\bn {\bf 2.2 Proposition}  {\it Let $X$ be a projective or compact K\"ahler manifold. Let $E$ be a numerically flat 
vector bundle on $X.$ Then there exists a filtration $$0 = E_0 \subset E_1  \subset \ldots \subset E_p = E$$ by vector subbundles
such that the graded pieces $E_{i+1}/E_i$ are unitary flat, i.e. defined by unitary representations of $\pi_1(X).$}

\bn Proposition 2.2 has the following consequence

\bn {\bf 2.3 Proposition}  {\it Let $E$ be a numerically flat vector bundle on the compact K\"ahler manifold $X.$
Let $h: \tilde X \la X $ be the universal cover. Then $h^*(E) = \O_{\tilde X}^r.$ } 

\bn {\bf Proof.} (a) The case that $E$ is unitary flat is obvious.
\sn (b) Next we consider the case of an exact sequence 
$$ 0 \la F \la E \la G \la 0 \eqno (*) $$
where both $F$ and $G$ are unitary flat. Then by (a) $h^*(F)$ and $h^*(G)$ are both trivial. Let $\zeta \in H^1(X,F \otimes G^*) $ be the 
extension class  defining (*). Then we must show that $h^*(\zeta) = 0.$ Since $F \otimes G^*$ is unitary flat, this
follows from the 

\bn {\bf 2.4 Fact} Let $V_r$ be a unitary flat bundle, $\zeta \in H^1(X,V).$  Then $h^*(\zeta) = 0.$
\sn (Notice however that $H^1(\tilde X,\O_{\tilde X}) $ might not vanish).
\bn {\bf Proof of the fact}: $V$ is given by a local system $\bE$. Hodge theory therefore yields an epimorphism
$$ H^1(X,\bE) \la H^1(X,E).$$
So we have an epimorphism
$$ h^*H^1(X,\bE) \la h^*H^1(X,E).$$
Since $h^*(\bE) = \tilde X \times \C^r$, it follows 
$$H^1(\tilde X,h^*(\bE)) = 0, $$
hence $h^*(\zeta) = 0$ in $H^1(\tilde X,h^*(E)) = 0,$ too.
\bn (c) Now we approach the general case in the proof of 2.3. We have an exact sequence 
$$ 0 \la E_i \la E_{i+1} \la E_{i+1}/E_i \la 0, \  i \geq 1,$$
and by (b) and induction resp. our assumption, both $h^*(E_i)$ and $h^*(E_{i+1}/E_i)$ are trivial. Denoting
$$\zeta \in H^1(X,(E_{i+1}/E_i)^* \otimes E_i)$$
the extension of the above        exact sequence, we need to show that $h^*(\zeta) = 0.$  This is done
by considering the induced sequence
$$ 0 \la (E_{i+1}/E_i)^* \otimes E_{i-1}Ê\la (E_{i+1}/E_i)^* \otimes E_i \la (E_{i+1}/E_i)^* \otimes E_i/E_{i-1}Ê\la 0,$$ 
taking $h^*$ and proceeding inductively using 2.4

\bn Before we can treat projective bundles, we need some preparations. The first one is classical and due to
Ehresmann:

\bn {\bf 2.5 Proposition} {\it Let $\varphi: X \la Y$ be a surjective fiber bundle of complex manifolds with typical fiber $F$  such that
$$ 0 \la \varphi^*(\Omega^1_Y) \la \Omega^1_X \la \Omega^1_{X/Y} \la 0 $$
splits. This exhibits a vector bundle $E \subset \Omega^1_X$ which is mapped isomorphically onto $\Omega^1_{X/Y}.$
We say that the splitting is integrable or that $E$ is integrable if 
$$dE \subset E \wedge \Omega^1_X.$$
Equivalently, $V = (\Omega^1_X/E)^*$ is integrable in the usual sense, i.e. closed under the Lie bracket. 
In that case there is a representation $\rho: \pi_1(Y) \la {\rm Aut}(F)$ so that $\pi_1(Y)$  acts on $F$ via $\rho.$
\sn Moreover $X$ is the fiber bundle over $Y$ with fiber $F$ given by $(\tilde Y \times F)/\pi_1(Y).$ Finally
the splitting $\tilde X = \tilde Y \times \tilde F$ is diagonal with respect to the above splitting. }

\bn {\bf Proof.} [Eh50], cp. [Be99]. 

\bn {\bf 2.6 Proposition} {\it Let $C$ be an irreducible reduced compact curve and $E$ a vector bundle of rank $r+1$ on $C. $
Let $X = \bP(E)$ with projection $p: X \la C.$ Suppose that the sequence
$$ 0 \la p^*(\Omega^1_C) \la \Omega^1_{\bP(E)}Ê\la \Omega^1_{\bP(E)/C}Ê\la 0$$
splits. Then $E' = E \otimes {{\rm det}E \over {r+1}}Ê$ is numerically flat. }

\bn {\bf Proof.}Ê(a) In case $C$ is smooth, we notice that the splitting is automatically integrable so that (2.5)
applies. Hence $X$ is defined by a representation $$\pi_1(C) \la {\rm Aut}(\bP_r) = \bP Gl(r+1,\C),$$ i.e. $E$ is 
projectively flat. Using the exact sequence
$$ 0 \la \mu_{r+1} \la Sl(r+1,\C) \la \P Gl(r+1,\C) \la 0, $$
and tensoring $E$ suitably, we find an \'etale cover $f:Ê\hat C \la C$ of degree $r$ such that $\hat E = f^*(E)$ is 
given by the lifted representation
$$ \hat \rho: \pi_1(\hat C) \la Sl(r+1,\C).$$
Hence $\hat E$ is flat, in particular $E$ is nef. Since also ${\rm det}E$ is flat, $\hat E$ and $E$ are both numerically
flat. 
\sn (b) Now suppose $C$ singular and let $\hat C \la C$ be the normalisation. Let  $\hat E = h^*(E)$ and put 
$\hat X = \bP(\hat E)$ with projection $\hat p: \hat X \la \hat C.$ It is sufficient to show that 
$$ 0 \la \hat p^*(\Omega^1_{\hat C}) \la \Omega^1_{\hat X} \la \Omega^1_{\hat X/\hat C} \la 0$$
splits; then by (a) $\hat E$, hence $E$, is numerically flat after normalising to $c_1 = 0.$ This splitting
now is obtained by an easy diagram chase from the lifted splitting on $X.$ 

\bn Now we are able to treat general $\bP_r-$bundles.

\bn {\bf 2.7 Theorem} {\it Let $X$ be a projective manifold with $\bP_r-$bundle structure $\varphi: X \la Y.$
Suppose that 
$$ 0 \la T_{X/Y} \la T_X \la \varphi^*(T_Y) \la 0 $$
resp.
$$ 0 \la \varphi^*(\Omega^1_Y) \la \Omega^1_X \la \Omega^1_{X/Y}Ê\la 0 \eqno (S)$$
splits. Then 
\sn (1) $\tilde X = \tilde Y \times \bP_r.$ 
\sn (2) If the splitting (S) is integrable, then the splitting $\tilde X = \tilde Y \times \bP_r$ is diagonal
with respect to (S), hence $X = (\tilde Y \times \bP_r)/\pi_1(Y),$ where 
$\pi_1(Y)$ acts on $\bP_r $ via a representation $\rho: \pi_1(Y) \la {\rm Aut}(\bP_r)$.
\sn (3) Suppose $X = \bP(V)$ with a vector bundle $V$ of rank $r.$ Then the $\Q-$bundle $V \otimes {{\rm det}V^* \over {r+1}}$
is numerically flat. }

\bn {\bf Proof.} (I) First we assume that $X = \bP(V).$ Let $V' = V \otimes {{\rm det}V \over {r+1}}.$ In order to 
prove (3), it is sufficient to show that the $\Q-$bundle $V'$ is nef. So let $C \subset X$ be an irreducible
curve. Let $X_C = \varphi^{-1}(C) = \bP(V\vert C).$ Then 
$$ \Omega^1_{X_C} \simeq \varphi^*(\Omega^1_C) \oplus \Omega^1_{X_C/C},$$
as one checks easily. So by (2.6), $V'\vert C$ is numerically flat, in particular nef. So $V'$ is  nef. This proves (3).
\sn Now $(r+1){\rm det}V'$ is a well-defined flat line bundle. Therefore we find a finite \'etale cover $g: \tilde Y \la Y$
such that $g^*(\det V')$ is Cartier, or, in other words, the determinant of $g^*(V)$ is divisible in ${\rm Pic}(Y)$  by $r+1.$ So we
may assume that $\det V$ itself is divisible by $r+1$ in ${\rm Pic}(Y).$
\sn Then $h^*(V') = \O_{\tilde X}$ by (2.3), in particular $\tilde X = \tilde Y \times \P_r.$ This proves (1).
\sn Finally (2) follows from (2.5), independently whether $X$ comes from a vector bundle or not.
\sn (II) It remains to prove (1) in the general case. Since $\tilde Y$ is simply connected, there exists a vector bundle
$\tilde V$ over $\tilde Y$ such that $\tilde X = \bP_r(\tilde V)$ (the obstruction lies in
$$ H^2(\tilde Y,\Z_{r+1}) \simeq H_1(\tilde Y,\Z_{r+1}) = 0 ).$$
Now consider the vector bundle
$$ E = \varphi_*(-K_{X/Y}).$$
By 2.6 $E\vert C$ is numerically flat for every curve $C \subset Y.$ Hence $E$ itself is numerically flat. 
By 2.3, $h_Y^*(E) = \O_{\tilde Y}^N.$ 
Now we conclude by Lemma 2.7.1.

\bn {\bf 2.7.1 Lemma}Ê{\it Let $Y$ be a complex manifold and $X$ a $\bP_r-$bundle over $Y$ with 
projection $\varphi: X \la Y.$
Suppose that $\varphi_*(-K_{X/Y}) = \O_Y^{N+1}.$ 
Then $X \simeq Y \times \bP_{r}.$}

\bn {\bf Proof.} Let $L = -K_{X/Y}.$
The canonical epimorphism 
$$ \varphi^*(\varphi_*(L) =  \la \O(k)$$
induces a map
$$ X = \bP(L) \la \bP(\varphi_*(L)) = \bP(\O_Y^{N+1}) = Y \times \bP_N.$$
Let $f: X \la \bP_N$ be the induced map. 
Then $f$ is given by $N+1$ sections in $L$ whose images trivialise $\varphi_*(L) = \O_Y^{N+1}.$ 
Since $c_1(L)^{r+1} = 0,$ we have $\dim f(X) \leq r.$ On the other hand, $f \vert \varphi^{-1}(y)$ is an 
embedding for all $y \in Y,$ so $\dim f(X) = r$ and actually $f(X) = \bP_r.$ 
Now the fibers of $f$ and $\varphi$ meet transversally, so that the induced map
$$ X \la Y \times \bP_r$$
is an isomorphism.

\bn We return to our original aim to investigate contractions $\varphi: X \la Y$ of a smooth projective threefold
$X$ with $T_X = L \oplus V.$

\bn {\bf 2.8 Theorem} {\it Suppose $\dim Y = 2.$ Then $\varphi$ is a $\P_1-$bundle and hence $\tilde X = \tilde Y \times
\P_1.$ Moreover either
\item {(a)} $L = \varphi^*(L')$ with a line bundle $L'$ on $Y$ and $T_Y = L' \oplus V'$ with a line
bundle $V'$. In this case $V$ is integrable and $\tilde X$ splits diagonally with respect to $T_X = L \oplus V.$ Or
\item {(b)} $L = T_{X/Y}.$ In this case $V$ is not necessarily integrable.}

\bn {\bf Proof.} Assume that $\varphi$ is not a $\P_1-$bundle. Then the general singular fiber is a reducible conic $C_1 + C_2$
($C_1 \ne C_2.$) We have $T_X \vert C_i = \O(2) \oplus \O \oplus \O(-1).$ Then (1.3) implies $T_{C_i}Ê= L \vert C_i$ contradicting
the fact that $C_1 $ and $C_2$ meet tranversally. So $\varphi$ is a $\P_1-$bundle (analytically).
\sn (a) Suppose that the canonical map $\gamma: T_{X/Y}Ê\la L$ vanishes. Then $T_{X/Y}Ê\subset V$ and (1.3) implies that
$L_C = \O_C$ for all fibers $C$ of $\varphi$ and that $V \vert C = \O_C(2) \oplus \O_C(-2).$ 
Hence $L = \varphi^*(L')$ for some line bundle $L'$ on $Y.$ In other words, the tangent map $$T\varphi \vert L: L \la L'$$
is an isomorphism. Then the splitting $T_X \la L$ yields a splitting map $T_Y \la L'$.
Writing $ T_Y = L' \oplus Q'$, it is clear that we have a sequence
$$0 \la T_{X/Y}Ê\la V \la \varphi^*(Q') \la 0.$$
By Beauville [Be99], the universal cover $\tilde Y$ splits as $\tilde Y = A_1 \times A_2$ with projections  $p_i: A_i \to Y$
and such that $T_{A_1}Ê= p^*(L')$ and $T_{A_2}Ê= p^*(Q').$ Then the local product structure of $\tilde X = X \times_{\tilde Y}Y$
(i.e. locally over $\tilde Y$) yields the integrability of $V$ and therefore by (2.5) $\tilde X$ splits diagonally with
respect to $T_X = L \oplus V.$ 
\sn (b) If $\gamma \ne 0, $ then $\gamma$ is injective, and $L_C = \O_C$, $V_C = \O \oplus \O(-2)$ for all fibers
$C$ of $\varphi.$ Hence $T_{X/Y}Ê\la V$ vanishes and $L = T_{X/Y}.$ 
For a non-integrable example, see [Dr99].

\bn {\bf 2.9 Theorem} {\it Suppose $\dim Y = 1.$ Then $\varphi$ is a $\P_2-$bundle, and $\tilde X$ splits diagonally
with respect to $T_X = L \oplus V.$}

\bn {\bf Proof.} (a) First we show that $\varphi$ is a $\P_2-$ bundle. 
It is sufficient to know that the general fiber of $\varphi$ is $\P_2$
[Mo82]. If not, then $F$ is $P_1 \times \P_1$ or some other del Pezzo surface, just called del Pezzo for
simplicity. If $F$ is del Pezzo, choose a $(-1)-$curve $l \subset  F.$ Then
$$ T_X \vert l = \O(2) \oplus \O \oplus \O(-1). $$
Then (1.3) yields $L \vert l = \O_l,$ hence $L$ is $\varphi-$trivial. It follows that the canonical map $\gamma: 
T_{X/Y}Ê\la L$ vanishes, hence $T_{X/Y}Ê\subset  V.$  This however says that $T_{X/Y}Ê= V$, and since all fibers
of $\varphi$ are reduced ([Mo82]), $\varphi$ must be a submersion. Now apply the Leray spectral
sequence to obtain
$$ H^2(X,\Z) = H^2(Y,\Z) \oplus H^0(Y,R^2\varphi_*(\Z)).$$
Hence $b_2(X) = b_2(F) + b_2(Y) \geq 2 + b_2(Y)$ contradicting the fact that $\varphi $ is the contraction
of an extremal ray, which implies that $b_2(X) = b_2(Y) + 1.$ 
\sn So $F$ cannot be del Pezzo.
\sn The case that $F = \P_1 \times \P_1$ is ruled out in the same way by choosing $l$ to be a ruling
line in such a way that $T_l \ne L_l.$ 
\sn (b) So $\varphi$ is a $\P_2-$bundle. Applying once again (1.3) we see that $L \vert l = \O$ and
$V \vert l = \O(2) \oplus \O(1)$ for every line $l$ in a fiber of $\varphi.$ Thus $L = \varphi^*(T_Y)$ and
$V = T_{X/Y}$ so that we can apply (2.5) and conclude. 

\bn In total we have proved

\bn {\bf 2.10 Theorem}  {\it Let $X$ be a smooth projective threefold with $\kappa (X) = - \infty.$ Suppose $T_X = L
\oplus V.$ Then $\tilde X$ splits. The splitting is diagonal with respect to $T_X = L \oplus V$ unless the
following holds: $X$ is the successive blow-up along smooth curves in a smooth projective threefold $Y$,
the splitting $T_X = L \oplus V$ induces canonically a splitting $T_Y = L' \oplus V'$ and there is a 
$\P_1-$bundle structure $\psi: Y \la Z$ such that $L' = T_{Y/Z}.$ }

\bn In order to get results in higher dimensions one would have to consider contractions $\varphi: X \la Y$
of fiber type, assuming $T_X = E_1 \oplus E_2.$   One expects that one of the $E_i$ is $\varphi-$trivial. 
We come back to this situation in the next section when we consider Fano manifolds. 

\bn \bn 

{\medium \centerline {3. Fano Manifolds}}
\bn
In this section we fix a Fano manifold $X$ of dimension $n$ with a splitting $T_X = \bigoplus_{i=1}^{k} E_i.$
Arguing inductively and having in mind that Fano manifolds are simply connected, we may assume $k  = 2.$ Let $r_i = {\rm rk}E_i$ and put $L_i = {\rm det} E_i.$ 
\bn
{\bf 3.1 Proposition} {\it We have $b_2(X) \geq 2.$}

\bn {\bf Proof.} Suppose $b_2(X) = 1$ and write ${\rm det}E_i = \O_X(m_i),$ where $\O_X(1)$ is the ample generator
of ${\rm Pic}(X) = \Z.$ By (0.4) we have $c_1(E_i)^t = 0$ for $t > r_i,$ therefore $m_i = 0$ for all $i.$ 
On the other hand $0 < c_1(X) = \sum m_i,$ contradiction. 

\bn {\bf 3.2 Remark} The same proof shows that the tangent bundle of a projective manifold $X$ with $K_X$ 
ample and $\rho(X) = 1$ cannot split. If $K_X \equiv 0,$ the last remaining case if $\rho(X) = 1,$ 
Beauville's decomposition theorem applies and also easily yields a contradiction. Hence the tangent bundle
of a projective manifold with $\rho(X) = 1$ never splits.

\bn The significance of (1.6) in our context is given by

\bn {\bf 3.3 Theorem}  {\it Let $X$ be a Fano manifold with $T_X = E_1 \oplus E_2.$ Suppose that (1.6) has a 
positive solution (e.g. every contraction of $X$ has a non-trivial smooth fiber). 
Then $X \simeq Z_1 \times Z_2$ and $T_{Z_i} = p_i^*(E_i).$ }

\bn {\bf Proof.} Let $\varphi: X \to Y$ be some contraction defined by the ray $R = \bR [l].$ Then $L_i \cdot l = 0$
for some $i$, say $i = 1,$ by (1.3) and our assumption. Since $-K_X\cdot l > 0$ and since $-K_X = L_1 \otimes L_2,$
we conclude that $L_2 \cdot l > 0.$ 
\sn It follows that $L_i \cdot l_j \geq 0$ for all extremal curves $l_j$ and $i = 1,2.$
Hence $L_i$ are both nef by the cone theorem. Now by (0.4) we obtain
$$ (-K_X)^n = {n \choose r_1} L_1^{r_1}Ê\cdot L_2^{r_2}.$$
In particular $$L_i^{r_i}Ê\ne 0, L_i^{r_i+1}Ê= 0. \eqno (*)$$
Since by the base point free theorem both $L_i$ are semi-ample and therefore define morphisms
$\psi_i: X \la Z_i$ to  normal projective varieties $Z_i.$ Both $\psi$ are contractions of extremal faces.
By (*) we have $\dim Z_i = r_i.$
\sn We claim that $\psi := \psi_1 \times \psi_2: X \la Z_1 \times Z_2$ is an isomorphism. 
Since $L_1 \cdot C = 0$ implies $L_2 \cdot C >  0$ for every extremal rational curve and vice versa, this actually holds for
every curve $C$,               hence  $\psi_2\vert \psi_1^{-1}(z)$ is always
finite and vice versa. Therefore $\psi$  is finite onto its image. Since $\dim (Z_1 \times Z_2) = r_1 + r_2 = n,$
$\psi $ is surjective and finite. 
Let $F$ be a general fiber of $\psi_1,$ then $F$ is Fano. Let $C \subset F$ be a smooth rational curve with $T_F \vert C$ ample;
then $T_X \vert C$ is nef. Since $L_1 \vert C = \O_C$ and since  $E_1  \vert C$ is nef as quotient of $T_X \vert C$,
it follows that $E_1 \vert C = \O_C^{r_1}.$ Therefore the composed map
$$ T_F \vert C \la T_X \vert C \la  E_1 \vert C$$
vanishes (remember $T_F \vert C$ ample) and hence, varying $C$, $T_F$ is a subbundle of $  E_2 \vert F.$ 
Since ${\rm rk}T_F = {\rm rk} E_2 $, we conclude $T_F = E_2 \vert F.$ 
\sn Hence $E_2$ is generically integrable, hence integrable. All fiber of $\varphi_1$ are clearly leaves
of $E_2.$ By symmetry, $E_1$ is integrable, too, and all fibers of $\varphi_2$ are leaves of $E_1.$ 
Thus $\varphi_1^{-1}(x)$ and $\varphi_2^{-1}(y)$ are smooth and transversal at every intersection point.
In particular $\varphi_1\vert \varphi_2^{-1}(y)$ is \'etale (and surjective) over $Z_1.$ But $\varphi_2^{-1}$ is Fano,
hence simply connected and so $\varphi_1 \vert \varphi_2^{-1}(y) $ is an isomorphism. This same holds for
$\varphi_2 \vert \varphi_1^{-1}(x)$ and thus $\psi $ is an isomorphism. The proof also shows that the splitting
is diagonal.

\bn We have used (1.6) only in order to make sure that the line bundles $L_i$ are $\varphi-$nef for all contractions $\varphi$,
i.e. that the $L_i$ are nef. Therefore we can state

\bn {\bf 3.4 Proposition}  {\it If $L_1 $ and $L_2$ are nef (in (3.3), then $X \simeq Z_1 \times Z_2$ and the splitting
is diagonal.}

\bn Of course the nefness of the $L_i$ is also necessary for the splitting.
\sn We next show that we have splitting (without assuming (1.6)) if the rank of say $E_1$ is small.

\bn {\bf 3.5 Theorem} {\it Let $X$ be Fano with $T_X = E_1 \oplus E_2$ If ${\rm rk}E_1 = 1,$ then $X \simeq Z_1 \times
Z_2$ diagonally.}

\bn {\bf Proof.} By (3.4) we only need to show that $L_1$ and $L_2$ are nef. Suppose that one $L_i$ is not $\varphi-$nef
for some contraction $\varphi.$ Since $-K_X = L_1 + L_2,$ it follows that $L_1$ is not $\varphi-$trivial. Hence $\dim \varphi^{-1}(y) \leq 1 $ for all $y$ by (1.8)(1). Let $F$ be an irreducible (reduced)
component of $\varphi^{-1}(y).$ Then $F \simeq \P_1$ since $R^1\varphi_*(\O_X) = 0.$ Moreover $\varphi^{-1}(y) = F,$
see [AW97,1.6]. Hence we can apply (1.8)(3) and derive a contradiction.

\bn {\bf 3.6 Theorem}Ê \ {\it Let $X$ be Fano with $T_X = E_1 \oplus E_2$. If ${\rm rk}E_1 = 2,$ then $X \simeq Z_1 \times Z_2$
diagonally. }

\bn {\bf Proof.} We proceed as in (3.5) and assume that some $L_i$ is not $\varphi-$nef. So $L_1$ is not $\varphi-$trivial. Let $F$ be a non-trivial fiber (with 
reduced structure). Then $\dim F \leq 2$ by (1.8). If $F$ contains a component of dimension 1, we conclude again by
[AW97,1.6]. Hence $F$ is purely 2-dimensional. By [AW97,1.19], any irreducible component $F_0$ of $F$ is either a
$\P_2$ or a Hirzebruch surface $\P(\O \oplus \O(-k)), $ or a Hirzebruch surface ($k \geq 2$) with the exceptional section blown down. 
In the first two 
cases we conclude by (1.8). The last requires some more work. Let 
$$\tilde F = \P(\O \oplus \O(-k)) \la F_0$$ be the minimal desingularisation.
Let $C_0 \subset  \tilde F$ be the exceptional section, $\tilde l \subset \tilde F$ a ruling line and set $l = \sigma (\tilde l).$
We claim that $l $ is a standard rational curve; then we will be done again by (1.8). 
\sn We consider the exact sequence
$$ 0 \la (T_F\vert l)/{\rm tor} {\buildrel {\alpha} \over {\la}} T_X \vert l \la N_{F\vert X}\vert l.$$
Since $N^*_{F/X}\vert l$ is generically spanned, we only have to investigate the inclusion 
$ V := (T_F \vert l)/{\rm tor}Ê\subset  T_X \vert l.$ 
Observe that the epimorphism
$$ \Omega^1_X \vert l \la \Omega^1_F \vert l \la 0$$
yields a vector bundle epimorphism
$$ \Omega^1_X \vert l \la (\Omega_F \vert l)/{\rm tor} \la 0.$$ 
Therefore $$V = (\Omega^1_F \vert l)/{\rm tor})^{*}Ê$$ is a subbundle of $T_X \vert l.$ 
So we obtain a vector bundle sequence
$$ 0 \la V \la T_X \vert l \la {\rm Coker}\alpha \la 0 $$
with $({\rm Coker}\alpha)^*$ being nef. Hence things come down to prove
$$ V \vert l = \O(2) \oplus \O(a)  \eqno (*)$$
with $a \leq 1.$ 
In order to verify (*) we consider the analytic preimage $\sigma^{-1}(l) = \tilde l + \mu(k) C_0.$ 
Then we have an epimorphism
$$ \sigma^*(N^*_{l/F}) \la N^*_{\sigma^{-1}(l)/ \tilde F}Ê\la 0;$$
hence 
$$ (\sigma^*(N^*_{l/F}/{\rm tor})Ê\vert \tilde l = N^*_{\sigma^{-1}(l)/\tilde F}Ê\vert \tilde l.$$
Now we use the exact sequence
$$ 0 \la N^*_{\sigma^{-1}(l)/\tilde F} \vert \tilde l \la N^*_{\tilde l/\tilde F}Ê\la N^*_{\tilde l/\sigma^{-1}(l)} \la 0, $$
and observe that $N^*_{\tilde l/\tilde F}Ê= \O.$ Since $N^*_{\tilde l/ \sigma^{-1}(l)}$ is supported at $\tilde l \cap C_0$, which
is just one point and since the stalk at that point has dimension 1, we conclude
$$ N^*_{\sigma^{-1}(l)/\tilde F}Ê\vert \tilde l = \O(-1).$$ 
Hence $\sigma^*(N^*_{l/F}/{\rm tor}) = \O(-1) $ by (*) and thus 
$$ N_{l/F}/{\rm tor} = \O(1).$$ 
From 
$$ 0 \la T_l \la T_F \vert l \la N_{l/F}Ê$$
we obtain an exact sequence
$$ 0 \la T_l \la  V = (T_F \vert l)/{\rm tor} \la N_{l/F}/{\rm tor}Ê= \O(1)$$
and thus $V = \O(2) \oplus \O(a)$ with $a \leq 1,$ proving (*) and therefore the theorem.     

\bn {\bf 3.7 Corollary} \ Ê{\it Let $X$ be a Fano $n-$fold, $n \leq 5.$ If $T_X = E_1 \oplus E_2,$
then $X \simeq Z_1 \times Z_2$ diagonally. }

\bn {\bf 3.8 Remark} \ Notice that the  proof of (3.6) provides a positive answer to the problem
1.6 in case $\varphi$ is birational with some non-trivial fiber of dimension at most 2. So in dimension 4
it only remains to treat the case that $\varphi$ contracts a divisor to a point.

\vfill \eject

{\medium \centerline {4. The Case $\kappa (X) = 3$}}
\bn 
In our study of threefolds $X$ with splitting tangent bundles the previous sections reduce us to $K_X$ nef. If 
$\kappa (X) = 0,$ then $mK_X = {\cal O}_X $ for a suitable positive $m$. Then by the decomposition theorem [Be83], we obtain
the diagonal splitting of $\tilde X$ we are looking for. This holds in all dimensions. In this section we treat the case that
$\kappa (X) = \dim X = 3,$ i.e. $K_X$ is 
big and nef in connection with our splitting
$T_X = L \oplus V.$ The remaining cases $\kappa (X) = 1,2$ will be treated in the second part of this paper. 
\sn By [Be99] (see 0.5) we may assume that  $K_X$ is not
ample (if $K_X$ is ample, then $X$ is K\"ahler-Einstein), therefore we have a non-trivial morphism (with connected fibers) $\phi: X \la Y$ to the canonical 
model of $X,$ given by the base point free linear system $\vert mK_X \vert $ for suitable large $m.$ 
The threefold $Y$ has only canonical singularities and the $\Q-$divisor $K_Y$ is ample.
We denote by $B \subset Y$ the exceptional locus of $\phi$ in $Y,$ this is nothing than the singularity
set of $Y.$ Furthermore let $E = \phi^{-1}(B).$
\bn Our aim is to prove

\bn {\bf 4.1 Theorem} {\it Let $X$ be a smooth projective threefold $K_X$ big and nef such that $T_X = L \oplus V$; 
let $\phi: X \la Y$ be the canonical model.
{\item (1) The universal cover $\tilde Y$ of $Y$ is of the form $\tilde Y \simeq \Delta \times S$ with $\Delta \subset \C$ 
the unit disc and $S$ a surface with only rational double points as singularities.
\item (2) If $g: \hat S \la S$ denotes the minimal resolution, then the universal cover $\tilde X$ of $X$ is of the form 
$$\tilde X \simeq \Delta \times
\hat S$$ and $\phi \simeq {\rm id}_{\Delta}Ê\times g.$  This decomposition is compatible with the decomposition
$T_X = L \oplus V,$ i.e. $\tilde X$ splits diagonally with respect to $T_X = L \oplus V.$ }}

\bn So locally (in $Y$) we have a product structure $Y \simeq \Delta \times W$ with a surface $W$ having only rational 
double points so that (semi-locally) $X \simeq C \times \hat W$ with $\hat W \la W$ the minimal resolution. 
In particular $B = C \times {\rm Sing}W$ with a smooth curve $C,$ i.e. $B$ consists of a disconnected
union of smooth curves. Our strategy will be to first establish this local product structure in $Y;$ in particular
this shows that $Y$ is an orbifold, i.e. $Y$ is locally of the form $\C^3 /G.$ Then we will use the theory
of orbifolds to conclude. 

\bn We start with an easy lemma which we will use frequently.

\bn {\bf 4.2 Lemma} {\it Let $X$ be a smooth complex threefold with $mK_X = \O_X$ for some positive integer $m.$ Let 
$Y$ be a normal Stein space and $f: X \la Y$ be a proper surjective bimeromorphic map. Suppose that the fibers of $f$ have dimension
at most $1.$ Suppose furthermore $R^1f_*(\O_X) = 0.$ Let $C \subset f^{-1}(y)$ be an irreducible curve. Then
$C \simeq \P_1$ with normal bundle $N_{C \vert X}Ê= \O(a) \oplus \O(b)$ and $(a,b) = (-1,-1),(0,-2),(1,-3).$}

\bn {\bf Proof.} We follow the arguments of Laufer [La81].  Let $\I$ denote the ideal of $C \subset X.$ Since $Y$ is Stein,
our vanishing assumption immediately gives
$$ H^1(X,\O_X) = H^1(Y,\O_Y) = 0.$$ 
Therefore $H^1(X,\I) = 0$ 
via $$ H^0(X,\O_X) \la H^0(C,\O_C) \la H^1(X,\I) \la 0.$$
Moreover an easy Leray spectral sequence argument yields
$$ H^2(X,\I^2) = 0$$ 
(here we use the assumption that there are no 2-dimensional fibers). 
Thus the exact sequence
$$ H^1(X,\I) \la H^1(C,\I/ \I^2) \la H^2(X,\I^2)$$
proves $$H^1(C,N_C^*) = 0. \eqno (*)$$
Let $F$ denote the reduced fiber $f^{-1}(y).$ Then the vanishing $R^1\phi_*(\O_X) = 0$ implies (via the theorem on formal functions and the
vanishing $H^2(F,\I^{\mu}Ê/ \I^{\mu+1})$ for all $\mu \geq 1$) that
$H^1(F,\O_F) = 0$ (having still in mind $\dim F = 1.$) 
Hence $H^1(C,\O_C) = 0$ and $C \simeq \P_1.$ Since $K_X \vert C = 0$ by our assumption on $K_X$, we have
$a + b = -2.$ Now (*) implies the claim on the normal bundle.

\bn {\bf 4.3 Proposition}  {\it $L$ and $V$ are integrable.}

\bn {\bf Proof.}  By [En87] $T_X$ is $K_X-$semistable. Then it is immediately verified that 
$$ c_1(L)\cdot K_X^2 = {1 \over {2}}Êc_1(V) \cdot K_X^2 = {1 \over {3}} K_X^3. $$
Therefore by [Su90,7.1] every $x \in X \setminus E$ has a 
neighborhood $U$ of the form $U = V \times V' $ such that $$T_U = p^*(T_V) \oplus p'^*(T_{V'}).$$
Hence $L$ and $V$ are generically integrable, hence everywhere integrable (the integrability of $L$ is anyway clear by reasons of rank). 

\bn {\bf 4.4 Corollary} {\it Every point $x \in X$ admits a neighborhood $U = U_1 \times U_2$ such that (by abuse of notation)
$T_{U_1}Ê= L$ and $T_{U_2}Ê= V.$}

\bn We now establish the local product structure in $Y.$

\bn {\bf 4.5 Proposition} {\it The exceptional locus $B \subset Y$ is smooth and for every $b \in B$ there are open neighborhoods $U$ of $b \in Y$ and
$V$ of $b \in B$ and there is a surface $W$ admitting a rational double point as only singular point such
that $U \simeq V \times W.$ In particular $Y$ is an orbifold.} 

\bn {\bf Proof.} {\bf (4.5.1)} In a first step we show that $\phi $ does not contract divisors to points. So assume to the
contrary the existence of an irreducible divisor $E_0 \subset E$ such that $\dim \phi(E_0) = 0.$ First we claim that
the canonical map $$ N_{E_0}^* \la L^*_{E_0} $$
obtained by composing $N^*_{E_0} \la \Omega^1_X \vert E_0$ and $\Omega^1_X \vert E_0 \la L^*_{E_0} $ vanishes.
In order to verify this claim, notice that
$$ H^0(N_{E_0}^{*k}) \ne 0$$
for suitable large $k$ ($E$ is contained in a fiber!). It follows easily that any desingularisation of $E_0$ has negative
Kodaira dimension; just observe that sections in $(N^*_E)^{k_0}        $ must have zeroes.
Now we choose a covering family $(l_t)$ of (generically irreducible) rational curves in $E_0.$ Then $E_0 \cdot L_t \leq 0.$
If $E_0 \cdot l_t = 0,$ we see easily that - denoting $f$ the normalisation of $l_t$ - $f^*(T_X)$ is nef, hence 
the deformations of $l_t$ cover $X$ which is absurd. Hence $$E_0 \cdot l_t < 0. \eqno (*)$$ 
Next notice $L \cdot l_t \geq 0.$ In fact, if $L \cdot l_t < 0,$ then for the normalisation $f_t: \P_1 \la l_t$ we obtain
$$H^1(\P_1,f_t^*(L^*_t)) = 0.$$ This implies $c_1(f_t^*(L_t)) = 0,$ contradiction. 
Combining $L_\cdot l_t \geq 0$ with (*), we conclude that the map $N^*_{E_0}Ê\to L^*_{E_0}$ vanishes. Hence $L_{E_0} \subset T_{E_0}.$
This means that at least generically $L_{E_0}$ defines a foliation on $E_0.$ 
\sn Now we apply (4.4) taking over all notations.
It follows that $ \tilde U = E_0 \cap U$ has a projection $\tilde U \la C$ to a possibly non-normal curve $C \subset U_2,$ and
moreover - possibly after shrinking - $\tilde U = U_1 \times C.$ 
We observe that
$$ \Omega^1_{E_0} = L^*_{E_0} \oplus W^*,$$
where $W^*$ is the image of $V^*_{E_0}Ê\la \Omega^1_{E_0}.$ Locally of course, $W^* = \Omega^1_C.$ 
 
Now we study the normalisation $h: \tilde E \la E_0$ of $E_0.$ Then $\tilde E $ is smooth since locally $\tilde E = U_1 \times
\tilde C,$ where $\tilde C$ is the normalisation of $C.$ 
If we introduce $\tilde L = h^*(L_{E_0}) $, then $\tilde L$ is an integrable subbundle
of $T_{\tilde E.}$ This is clear from the local description of $E.$ Let $Q$ denote the quotient bundle. 
Then generically $Q$ is a subbundle of $T_{\tilde E}.$ Since every point $x 
\in \tilde E$ has a neighborhood of the form $U_1 \times \tilde C$ (shrinking $\tilde C$, of course), such that $T_{U_1} = \tilde L,$
we conclude that $T_{\tilde C}Ê= Q$ and therefore $Q$ is a subbundle of $T_{\tilde E}$ everywhere, hence
$$ T_{\tilde E}Ê= \tilde L \oplus Q.$$ Since $\kappa (\tilde E) = - \infty,$ we thus have a $\P_1-$bundle structure $\tilde p: \tilde E \la \tilde B$
over a smooth curve $\tilde B$ and either $L = T_{\tilde B}$, $Q = T_{\tilde E / \tilde B}$ or vice versa [Be99]. 
Moreover we can write $\tilde E = \P(\E) $ with a numerically flat 
rank 2 bundle $\E$ on $\tilde B$ (2.6). Let $C_0 \subset  \tilde E$ be a section with $C_0^2 $ minimal. Since $\E$ is flat,
we have $C_0^2 = 0,$ cp. [Ha77,V.2]. Let $\tilde F$ be a fiber of $\tilde p.$ 

\bn (A)  Suppose first that $E$ is smooth, hence $\tilde E = E.$ Using the canonical vector bundle sequence 
$$ 0 \la N^*_{E_0} \la V^* \vert E_0 \la W^* \la 0$$
and $$ K_X \vert E_0 = L^* \vert E_0 + \det V^* \vert E_0 = 0,$$
we obtain $W^* = K_{E_0 / \tilde B}$ in case $L = T_{\tilde B}$ and $W^* = p^*K_{\tilde B}$ in case $Q = T_{\tilde B}.$
In both cases we clearly have $H^2(E_0,\det V^*) = 0 $ and since
$$ c_2(V \vert E_0 ) \in H^2(E_0,\det V^*),$$
we conclude $c_2(V \vert E_0 ) = 0$. On the other hand, $c_2(V \vert E_0)$ can be computed by the above
exact sequence and we easily obtain $g(\tilde B) = 1.$ 

\bn (B) Now we treat the case $\tilde L = T_{\tilde B}.$ Write
$$ h^*(K_{E_0}) = aC_0 + b \tilde F. $$
Since $K_{E_0} = N_{E_0} $ and since $E \cdot F < 0$ for $F = h(\tilde F),$ we have $a < 0.$ 
Let $\tilde N \subset  \tilde E$ be the preimage of the non-normal locus of $E.$ Write
$$ \tilde N = cC_0 +d\tilde F.$$
Then $c \geq 0$ and either $\tilde N = \emptyset $ or $cÊ\geq 2$ (this is immediate by considering the subadjunction formula of a 
Gorenstein curve). 
Since $K_{\tilde E} \equiv -2C_0 + (2g-2)\tilde F,$
where $g = g(\tilde B),$ the subadjunction formula
$$ K_{\tilde E} = h^*(K_{E_0}) - \tilde N$$
yields $c \leq 1,$ hence $c = 0$ and $E$ is normal, hence smooth. So $\tilde E = E_0$ with a map $p: E_0 \la \tilde B$ to
the elliptic curve $\tilde B$ (by (A)). Now we contract $X$ along $p: E_0 \la \tilde B$ and obtain a birational
morphism $$\psi: X \la X'$$
with $\psi(E_0) = \tilde B.$ 
Note that $X'$ has canonical singularities, that $K_{X'}$ is big and nef and that there is a birational map $\tau: X' \la Y$
such that $\phi = \tau \circ \psi.$ In particular
$$ R^q\tau_*(\O_{X'}) = 0 \eqno (+) $$ for $q \geq 1. $
We claim that,  $\tilde B$ being elliptic, $\tilde B$ cannot be an irreducible component of the fiber $F' = \tau^{-1}(\tilde B) = \tau^{-1}(y).$  
In fact, suppose the contrary.
Let $F_0$ be the closure of $F' \setminus \tilde B$ in $F'.$ Let $\hat F'$ denote completion of $X'$ along $F',$ analogously 
$\hat B$ and $\hat F_0.$ Let $\hat \I$ be the ideal sheaf of $\hat B \subset  \hat F$ and $\hat \I_0$ that one of
$ \hat { (B\cap F_0  )} \subset \hat F_0.$ Then we immediately see that
$$ H^2(\hat F,\hat \I) = H^2(\hat F_0,\hat \I_0) = H^2(\hat F_0,\O_{\hat F_0}).$$
The last group vanishes because $R^2\tau_*(\O_{X'}) = 0.$ 
Via the sequence
$$ H^1(\hat F,\O_{\hat F}) \la H^1(\hat B,\O_{\hat B}) \la H^2(\hat F,\hat \I) $$
and via the obvious epimorphism
$$ H^1(\hat B,\O_{\hat B}) \la H^1(\tilde B,\O_{\tilde B}),$$
we obtain $H^1(\hat F,\O_{\hat F}) \ne 0,$ contradicting (+).
\sn Thus we obtain an irreducible divisor $E'' \subset F'$ containing $\tilde B.$ Let $E' \subset X$ be the strict 
transform in $X.$ Consider an irreducible curve $C \subset E_0 \cap E'.$ Since we also have $L \vert E' \subset T_{E'},$
the $L-$leaf through any $x \in C$ is contained in $E_0 \cap E'.$ Therefore $C$ must be a leaf of $L$ and $E_0 \cap E'$ has to
consist of disjoint sections of $E_0 \la \tilde B.$ Moreover $E'$ has the same structure as $E_0$ and we have a map
$p': E' \la \tilde B$ with the $C_i$ also being section of $p'.$ We consider the map
$$ q: E_0 \cup E' \la \tilde B.$$ 
It is easily verified that we can blow down $X$ along $q$ to obtain a new birational map $\psi': X \la X''.$
E.g. blow down $X'$ along $E'' \la \tilde B$ and compose with $\phi.$ Arguing with $\psi'$ as with $\psi$ we deduce 
the existence of a third divisor in $\phi^{-1}(y).$ This procedure must terminate, so that finally we find  birational maps
$$\lambda: X \la Z, \ \sigma: Z \la Y$$ 
with $\sigma \circ \lambda = \phi$ such that $\sigma^{-1}(y)$ contains an elliptic curve as irreducible component, which 
is a contradiction as before with $\tau$.
This settles the case $\tilde L = T_{\tilde B}.$

\bn (C) Now suppose $Q = T_{\tilde B}.$ Here $\tilde p: \tilde E \la \tilde B$ descends to a $\bP_1-$bundle $p: E_0 \la B_0$
with $\tilde B \la B_0$ the normalisation and $\tilde E = E \times_{B_0}Ê\tilde B.$ 
\sn Let $\psi: X \la X'$ be the contraction
along $p.$ Arguing as in (B) either $B_0 \simeq \P_1$ or there exists a divisor $E' \subset \phi^{-y}$ (where $y = \phi(E_0)$)
with $ B_0 = \psi(E_0) \subset\psi(E').$ Then $E_0 \cap E'$ contains some multi-section $C$ of $p.$ Take $x \in C$ and
consider $F_x = p^{-1}(x).$ Then $F_x$ is the $L-$leaf through $x.$ Again by $ L \vert E' \subset T_{E'},$ we must have
$F_x \subset E'.$ This is absurd for general $x$, hence $B_0 = \P_1.$ In particular $E$ is smooth and we conclude by (A).

\bn {\bf (4.5.2)} Now we know that all non-trivial fibers of $\phi$ are 1-dimensional. By (4.2), all components $C$
of these fibers are smooth rational curves. Moreover 
$$ T_X \vert C = \O(2) \oplus \O(a) \oplus \O(b) $$
with $(a,b) = (-1,-1), (0,-2), (1,-3).$ From (1.3), it follows that $L_C = \O(2)$ or $L_C = \O_C.$ 
The case
$L_C = \O_C(2),$ however, is ruled out again by 1.3. Thus $L_C = \O_C$ (so $a = 0$ and $b = -2$) and $T_C \subset V_C$, hence 
$C$ is contained in a leaf $\V_0$ of $V.$ Now apply (4.4): every $x \in C$ has an open neighborhood  $U(x) = U_1(x) \times
U_2(x)$ with $T_{U_1} = L$ and $T_{U_2}Ê= V.$ 
By conpactness of $C$ we obtain $x_1, \ldots, x_r$ such that
$$ C \subset \bigcup_i U(x_i) =: U.$$
Let $U_2 := \bigcup_i U_2(x_i)$; then $U$ is an open neighborhood of $C \subset \V_0.$ We obtain a projection
$$Ê\pi : U \la U_2$$
by setting $\pi (x) := \L_x \cap \V_0$, where $\L_x$ is the $L-$leaf containing  $x \in U.$ The fibers of $\pi$ are open
parts of $L-$leaves. Fixing a point $x_0 \in C,$ we clearly may assume that every leaf $\V$ with $\V \cap U = 
\not \emptyset$ meets $\L_{x_0}.$ Hence we find an open neighborhood $W$ of $x_0$ in $\L_{x_0}$ and a holomorphic map
$ \tau : U \la W$ by setting $\tau (x) = \V_x \cap \L_{x_0},$ $\V_x$ denoting the $V-$leaf through $x_0.$   
After a coordinate change, $\tau^{-1}(0) = \V_0 \cap U$ and $W = \Delta_{\epsilon}Ê= \{z \in \CÊ\vert \  \vert z \vert < \epsilon\}.$
Since the fibers of $\tau$ and $\pi$ meet only tranversally (the leaves $\L_x$ and $\V_y$ meet only transversally), 
the map $$\tau \times \pi: U \la \Delta_{\epsilon}Ê\times \tilde U$$ is biholomorphic. 
Therefore $C$ obviously deforms to the leaves $\tau^{-1}(t)$. Hence $E$ cannot contain any 1-dimensional irreducible
component, so $\phi$ contracts only divisors to curves.
If $E_0 \subset E$ is an irreducible component, then we deduce that $E_0 \la \phi(E_0)$ is a $\P_1-$bundle and every 
fiber is a $(-2)-$curve some leaf $\V.$ To get the complete picture, consider the whole fiber $\phi^{-1}(y)$ for
$y \in B.$ Write
$$ \phi^{-1}(y) = C_1 + \ldots + C_r.$$
Then all arguments made above for a single $C$ also apply to $C_1 + \ldots + C_r$ as well. Identifying $U$ with 
$\tilde U \times \Delta_{\epsilon},$ it follows that $\phi\vert U = \psi \times {\rm id},$ where $\psi$ is the blow-down
of $C_1 + \ldots + C_r$ in the leaf $\V_0$. This completes the proof of (4.5).

\bn {\bf Proof of Theorem 4.1} By (4.4) $Y$ is an orbifold. Hence by Tian-Yau's result [TY87,] Theorem 2.2, there exists a
K\"ahler-Einstein-orbifold metric $\omega$ on $Y;$ so $[\omega]Ê\in c_1(Y).$ 
\sn Let $L' = \phi_*(L)^{**}$ and $V' = \phi_*(V)^{**}.$ Then
$$ T_YÊ= L' \oplus V'.$$
Now this decomposition is orthogonal with respect to $\omega.$ This is seen as in the case of an ordinary manifold, following
the arguments e.g. in [Ko87,p.178/179] (proof of theorem V.8.3). The crucial vanishing result [Ko87,III.1.9] immediately 
generalises to the orbifold case. 
Finally the de Rham decomposition for orbifolds (0.7)
applies and yields Part (1) of (4.1). Part (2) is then an obvious consequence of (1).

\vfill \eject 

{\centerline {\medium References}}
\bn

{Ê \item {[AW97]}ÊAndreatta,M;Wisniewski,J.:  A view on contractions of higher dimensional varieties.
Proc. Symp. Pure Math. 62(pt.1), 153-183 (1997) 

\item {[Be99]}ÊBeauville,A.: Complex manifolds with split tangent bundles. Complex analysis and algebraic geometry, vol. in memory
of M.Schneider, 61-70. de Gruyter 2000

\item {[DPS94]}ÊDemailly,J.-P.; Peternell,Th.; Schneider,M.:Compact complex manifolds with numerically effective tangent
bundles. J.Alg.Geom.3, 295--345 (1994)

item {[DPS00]} Demailly,J.-P.; Peternell,Th.; Schneider,M.: Pseudo-effective line bundles on compact K"ahler manifolds. 
math.AG/0006205 , to appear in Intl. J. Math. 

\item {[Dr99]}ÊDruel,S.: Vari\'et\'es alg\'ebriques dont le fibtr\'e tangent est totalement decompos\'e. J. reine u. angew.
Mathematik 522 (2000) 

\item {[Eh50]}ÊEhresmann,C.: Les connexions infinit\'esimales dans un espace fibr\'e diff\'erentiable. 
Colloque de topologie, Bruxelles (1950). 29-55. G.Thone, Liege 1951

\item {[En87]}ÊEnoki,I.: Stability and negativity for tangent bundles of minimal K\"ahler spaces. Lecture Notes in Math.
vol. 1339, 118-127. Springer 1987

\item {[Ha77]}ÊHartshorne,R.: Algebraic Geometry. Springer 1977

\item {[HW98]}ÊHwang,J.M.;Mok,N.: Rigidity of irreducible hermitian symmetric spaces of the compact type under
K\"ahler deformation. Inv. Math. 131, 393-418 (1998)

\item {[KN63]}ÊKobayashi,S.;Nomizu,K.: Foundations of differential geometry, vol.I,II. Wiley 1963

\item {[Ko87]} Kobayashi,S.: Differential geometry of complex vector bundles. Iwanami Shoten and Princeton Univ. Press,
1987

\item {[Ko96]}ÊKollar.J.: Rational curves on algebraic varieties. Springer 1996

\item {[La81]}ÊLaufer, H.:  On C$\bP_1$ as an exceptional set.
Recent developments in several complex variables, Proc. Conf., Princeton Univ. 1979, Ann. Math. Stud. 100, 261-275 (1981).

\item {[Mo82]} Mori,S.: Threefolds whose canonical bundles are not numerically effective. 
Ann. Math. 116, 133 - 176 (1982) 

\item {[Su90]} Sugiyama,K.: Einstein-K\"ahler metrics on minimal varieties of general type and an inequality between Chern numbers. 
Recent topics in differential and analytic geometry, Adv. Stud. Pure Math. 18-I, 417-443 (1990)

\item {[TY87]}ÊTian,G.;Yau,S.T.: Existence of K\"ahler-Einstein metric on complete K\"ahler manifolds and their
applications to algebraic geometry. Math. Aspects of String Theory, 574-628, World Scientific 1987

}

\vskip .5cm
\vskip .5cm
$$\matrix { {\rm  Frededric \  Campana} & & & & &  {\rm Thomas \  Peternell} \cr
{\rm Departement \  de  \ mathematiques} & & & & & {\rm Mathematisches \  Institut} \cr
{\rm Universite \  de \ Nancy} & & & & & {\rm Universitaet \  Bayreuth} \cr
{\rm BP 239} & & & & & & \cr
{\rm 54506 \  Vandoeuvre \  les \  Nancy} & & & & & {\rm 95440 \  Bayreuth} \cr
{\rm France} & & & & & {\rm Germany} \cr}$$

\end